\documentclass[showpacs,preprintnumbers,amsmath,amssymb]{revtex4}

\usepackage{graphicx}
\usepackage{bm}
\usepackage{epstopdf}
\usepackage{float}
\usepackage{dcolumn}
\usepackage{bm}
\usepackage{mathrsfs}
\usepackage{amsmath}

\begin{document}

\title{New indefinite integrals of Heun functions}
\author{D. Batic}
\email{davide.batic@ku.ac.ae}
\affiliation{%
Department of Mathematics,\\  Khalifa University of Science and Technology,\\ 
Abu Dhabi, United Arab Emirates 
} 
\author{O. Forrest}
\email{omar_forrest@yahoo.com}
\affiliation{%
Department of Mathematics,\\  University of the West Indies,\\
Kingston 6, Jamaica 
}
\author{M. Nowakowski}
\email{mnowakos@uniandes.edu.co}
\affiliation{
Departamento de Fisica,\\ Universidad de los Andes,\\ 
Cra.1E No.18A-10, Bogota, Colombia
}%
\begin{abstract}
We present a conspicuous number of indefinite integrals involving Heun functions and their products obtained by means of the Lagrangian formulation of a general homogeneous linear ordinary differential equation. As a by-product we also derive new indefinite integrals involving the Gauss hypergeometric function and products of hypergeometric functions with elliptic functions of the first kind. All integrals we obtained cannot be computed using Maple and Mathematica.
\end{abstract}

\maketitle

\section{Introduction}

Given the ordinary second-order linear differential equation
\begin{equation}\label{0}
y^{''}(x)+P(x)y^{'}(x)+Q(x)y(x)=0,
\end{equation} 
\cite{con1,con2} derived the indefinite integral
\begin{equation}\label{1}
\int\! f(x)\left[h^{''}(x)+P(x)h^{'}(x)+Q(x)h(x)\right]y(x)\,\mathrm{d}x=f(x)W(y,h)(x)+c,
\end{equation}
where a prime denotes differentiation with respect to the independent variable and
\begin{equation}\label{f}
f(x):=e^{\int\! P(x)\,\mathrm{d}x},
\end{equation}
$h$, $P$, and $Q$ are complex-valued differentiable functions in the variable $x\in\mathbb{R}$ with $h$ at least twice continuously differentiable, and $W(y,h)(x)=y(x)h^{'}(x)-h(x)y^{'}(x)$ denotes the Wronskian. Possible methods allowing to obtain an indefinite integral involving the solution of (\ref{0}) by using (\ref{1}) are
\begin{enumerate}
\item
a choice of $h$ in terms of a simple elementary function such as 
\begin{equation}\label{scelta_h}
h(x)=x^m e^{\rho x^\ell}\left\{\begin{array}{c}
\sin{(k_1 x)}\\
\cos{(k_2 x)}
\end{array}
\right\}
\end{equation}
with $m,\ell\in\mathbb{N}_0:=\mathbb{N}\cup\{0\}$, $\rho,k_1,k_2\in\mathbb{C}$ and so on. Note that the above choices contain the more simple cases $h(x)=1$ and $h(x)=x^m$. It is worth to mention that when $h$ is a constant function, (\ref{1}) takes the simpler form
\begin{equation}\label{zero}
\int\! f(x)Q(x)y(x)\,\mathrm{d}x=-f(x)y^{'}(x)+c.
\end{equation}
\item
If we want to integrate a solution $y(x)$ of (\ref{0}), we need to choose $h$ so that
\begin{equation}\label{y_alone}
h^{''}(x)+P(x)h^{'}(x)+Q(x)h(x)=\frac{1}{f(x)}.
\end{equation}
\item
A specification of $h$ as a solution of the differential equation (\ref{0}) with one or two  terms in (\ref{0}) deleted. For instance, we might require that $h$ satisfies one of the following differential equations
\begin{eqnarray}
h^{''}(x)+P(x)h^{'}(x)&=&0,\label{2}\\
P(x)h^{'}(x)+Q(x)h(x)&=&0,\label{3}\\
h^{''}(x)+Q(x)h(x)&=&0.\label{4}
\end{eqnarray}
In the case that $P$ and $Q$ consists of multiple terms we might also try to specify $h$ with $P$ or $Q$ with some of their subterms removed. 
\item
Take $h$ to be the solution of the equation conjugate to (\ref{0}), that is $h^{''}(x)+P(x)h^{'}(x)+\overline{Q}(x)h(x)=0$. This equation has the same $P$ as in (\ref{0}) but different $Q$. Since we can always construct a transformation of the dependent variable to make any two differential equations conjugate, this method allows to construct indefinite integrals of products of the solutions of the two conjugate equations according to the formula \cite{con1,con2}
\begin{equation}\label{intcon}
\int\! f(x)\left[Q(x)-\overline{Q}(x)\right]h(x)y(x)\,\mathrm{d}x=f(x)\left[h^{'}(x)y(x)-h(x)y^{'}(x)\right]+c.
\end{equation}
\end{enumerate}
Here, we derive new indefinite integrals involving solutions of the Heun equation \cite{ronv} which is represented by (\ref{0}) with
\begin{equation}\label{pq}
P(x)=\frac{\gamma}{x}+\frac{\delta}{x-1}+\frac{\epsilon}{x-a},\quad
Q(x)=\frac{\alpha\beta x-q}{x(x-1)(x-a)}
\end{equation}
where $\alpha,\beta,\gamma,\delta,\epsilon\in\mathbb{C}$ satisfy the Fuchsian condition $\epsilon=\alpha+\beta+1-\gamma-\delta$,  $a\in\mathbb{R}\backslash\{0,1\}$, and $q\in\mathbb{C}$ is the so-called accessory parameter.

\section{Construction of indefinite integrals}
Since the function $h$ appearing in (\ref{1}) is not fixed, we have in practice an unlimited number of cases for any given special function. In the following we treat only those choices of $h$ such that equation (\ref{1}) allows to derive new and interesting indefinite integrals. Last but not least, we checked that Maple and Mathematica software packages are unable to evaluate the indefinite integrals we computed in this paper. As a by-product we also obtain new indefinite integrals involving the Gauss hypergeometric function. These integrals do not appear to be listed in \cite{Grad}.  As in \cite{HR} we restrict the local solutions to the analytic ones around the origin, here denoted by $H_l(a,q;\alpha,\beta,\gamma,\delta;x)$ but extensions are possible for other solutions around other singularities. Let us choose $h$ according to (\ref{scelta_h}). Moreover, (\ref{f}) gives $f(x)=x^\gamma (x-1)^\delta(x-a)^\epsilon$. Hence, from (\ref{1}) we obtain
\[
\int\!x^{\gamma+m-2}(x-1)^{\delta-1}(x-a)^{\epsilon-1}e^{\rho x^\ell}
F(x)
H_l(a,q;\alpha,\beta,\gamma,\delta;x)\,\mathrm{d}x=
\]
\begin{equation}\label{Heun2}
x^{\gamma+m-1}(x-1)^{\delta}(x-a)^{\epsilon}e^{\rho x^\ell}
\left\{
\begin{array}{c}
\mathfrak{q}(x)\sin{(k_1 x)}+k_1 x\cos{(k_1 x)}H_l(a,q;\alpha,\beta,\gamma,\delta;x)\\
\mathfrak{q}(x)\cos{(k_2 x)}-k_2 x\sin{(k_2 x)}H_l(a,q;\alpha,\beta,\gamma,\delta;x)
\end{array}
\right\}+c
\end{equation}
with 
\begin{eqnarray}
F(x)&=&\left\{
\begin{array}{c}
x\mathfrak{p}_1(x,k_1)\cos{(k_1 x)}+\mathfrak{p}_2(x,k_1)\sin{(k_1 x)}\\
\mathfrak{p}_2(x,k_2)\cos{(k_2 x)}-x\mathfrak{p}_1(x,k_2)\sin{(k_2 x)}
\end{array}
\right\},\\
\mathfrak{q}(x)&=&(m+\rho\ell x^\ell)H_l(a,q;\alpha,\beta,\gamma,\delta;x)-xH^{'}_l(a,q;\alpha,\beta,\gamma,\delta;x),\\
\mathfrak{p}_1(x,k)&=&k\sum_{i=0}^2\mathfrak{a}_i x^i+2k\rho\ell x^\ell(x-1)(x-a),\\
\mathfrak{p}_2(x,k)&=&\sum_{i=0}^4\mathfrak{b}_i x^i+\rho\ell x^\ell\left[\sum_{i=0}^2\mathfrak{c}_i x^i+\rho\ell x^\ell(x-1)(x-a)\right]
\end{eqnarray}
and
\begin{eqnarray}
\mathfrak{a_2}&=&\alpha+\beta+2m+1,~ \mathfrak{a}_1=\mathfrak{a}_2+\delta(1-a)-\mathfrak{a}_0,~\mathfrak{a}_0=a(\gamma+2m),\\
\mathfrak{b}_4&=&-k^2,\quad\mathfrak{b}_3=k^2(a+1),~
\mathfrak{b}_2=-ak^2+\alpha\beta+m(\alpha+\beta+m),\\ \mathfrak{b}_1&=&m[\delta(1-a)-\alpha-\beta-m]-\mathfrak{b}_0-q,~\mathfrak{b}_0=am(m+\gamma-1),\\
\mathfrak{c}_2&=&\alpha+\beta+\ell+2m,\quad
\mathfrak{c}_1=\delta+a(1-\gamma-\delta+\ell+2m)-\mathfrak{c}_2,\\
\mathfrak{c}_0&=a&(\ell+\gamma+2m-1).
\end{eqnarray}
From (\ref{Heun2}) we can readily obtain new indefinite integrals for the Gauss hypergeometric function. To this purpose, we only need to observe that choosing appropriately the parameters of the Heun equation as in \cite{HR,Maier} yields 
\begin{eqnarray}
H_l(2,\alpha\beta;\alpha,\beta,\gamma,\alpha+\beta-2\gamma+1;x)&=&{}_{2}F_{1}\left(\frac{\alpha}{2},\frac{\beta}{2};\gamma;\mathfrak{h}(x)\right),\label{L1}\\
H_l\left(4,\alpha\beta;\alpha,\beta,\frac{1}{2},\frac{2}{3}(\alpha+\beta);x\right)&=&{}_{2}F_{1}\left(\frac{\alpha}{3},\frac{\beta}{3};\frac{1}{2};\mathfrak{f}(x)\right),\label{L2}\\
H_l\left(2,\alpha\beta;\alpha,\beta,\frac{\alpha+\beta+2}{4},\frac{\alpha+\beta}{2};x\right)&=&{}_{2}F_{1}\left(\frac{\alpha}{4},\frac{\beta}{4};\frac{\alpha+\beta+2}{4};\mathfrak{g}(x)\right)\label{L3}
\end{eqnarray}
with
\begin{equation}\label{hfg}
\mathfrak{h}(x)=x(2-x),\quad
\mathfrak{f}(x)=\frac{x}{4}(x-3)^2,\quad
\mathfrak{g}(x)=-4x(x-1)^2(x-2).
\end{equation}
Furthermore, 15.2.1 in \cite{abra} gives
\begin{eqnarray}
&&\frac{H^{'}_l(2,\alpha\beta;\alpha,\beta,\gamma,\alpha+\beta-2\gamma+1;x)}{t(x)}={}_{2}F_{1}\left(\frac{\alpha}{2}+1,\frac{\beta}{2}+1;\gamma+1;\mathfrak{h}(x)\right),\label{L4}\\
&&\frac{H^{'}_l\left(4,\alpha\beta;\alpha,\beta,\frac{1}{2},\frac{2}{3}(\alpha+\beta);x\right)}{s(x)}={}_{2}F_{1}\left(\frac{\alpha}{3}+1,\frac{\beta}{3}+1;\frac{3}{2};\mathfrak{f}(x)\right),\label{L5}\\
&&\frac{H^{'}_l\left(2,\alpha\beta;\alpha,\beta,\frac{\alpha+\beta+2}{4},\frac{\alpha+\beta}{2};x\right)}{r(x)}=
{}_{2}F_{1}\left(\frac{\alpha}{4}+1,\frac{\beta}{4}+1;\frac{\alpha+\beta+6}{4};\mathfrak{g}(x)\right)\label{L6}
\end{eqnarray}
with $t(x)=\alpha\beta(1-x)/(2\gamma)$, $s(x)=\gamma t(x)(3-x)/3$, $r(x)=\kappa(1-x)(2x^2-4x+1)$, and $\kappa=2\alpha\beta/(\alpha+\beta+2)$. Moreover, by means of (\ref{zero}) we obtain immediately the following result
\[
\int\!x^{\gamma-1}(x-1)^{\delta-1}(x-a)^{\epsilon-1}(\alpha\beta x-q)H_\ell\left(a,q;\alpha,\beta,\gamma,\delta;x\right)\,\mathrm{d}x=
\]
\begin{equation}\label{F12}
-x^{\gamma}(x-1)^{\delta}(x-a)^{\epsilon}H^{'}_\ell\left(a,q;\alpha,\beta,\gamma,\delta;x\right)+c.
\end{equation}
As a verification of the correctness of (\ref{F12}), we observe that if we define  $z=\mathfrak{h}(x)$ from which $x=1+\sqrt{1-z}$, and let $a=\alpha/2$, $b=\beta/2$, and $c=\gamma$, then (\ref{F12}) together with (\ref{L1}) and (\ref{L4}) reproduces 1.15.3.9 in \cite{Prud}, namely
\begin{equation}\label{PrudF}
\int\!z^{c-1}(1-z)^{a+b-c}{}_{2}F_{1}\left(a,b;c;z\right)\,\mathrm{d}z=
\frac{z^c}{c}(1-z)^{a+b-c+1}{}_{2}F_{1}\left(a+1,b+1;c+1;z\right)+\widetilde{c}.
\end{equation}
As a further independent check, if we let $z=\mathfrak{f}(x)$ from which
\[
x=2+\frac{1}{g(z)}+g(z),\quad
g(z)=\sqrt[3]{2\sqrt{z^2-z}+2z-1},
\]
and define $a=\alpha/3$, $b=\beta/3$, and $c=1/2$, then (\ref{F12}) together with (\ref{L2}) and (\ref{L5}) reproduces again (\ref{PrudF}). We also arrive at the same conclusion if we consider $z=\mathfrak{g}(x)$ with $x=1+(\sqrt{2+2\sqrt{1-z}}/2)$, make the identification $a=\alpha/4$, $b=\beta/4$, and $c=a+b+1/2$, and then, use (\ref{L3}) and (\ref{L6}). An integral involving $H_l(a,q;\alpha,\beta,\gamma,\delta;x)$ alone can be constructed by making the choice $h(x)=\alpha^{-1}(1-\alpha)^{-1}$ with $\alpha\neq 0,1$, $q=\alpha(1-\alpha)$, $\gamma=\epsilon=1$ and $\delta=0$. In this case, we find
\begin{equation}\label{HFe}
\int\! H_l(a,\alpha-\alpha^2;\alpha,1-\alpha,1,0;x)\,\mathrm{d}x=\frac{x(a-x)}{\alpha(1-\alpha)}H^{'}_l(a,\alpha-\alpha^2;\alpha,1-\alpha,1,0;x)+c.
\end{equation}
Note that even in this specialized case Maple/Mathematica is not able to solve the above indefinite integral. As a verification of the correctness of our method, we let $a=\alpha$, $b=1-\alpha$, and $c=1$, then (\ref{HFe}) together with (\ref{L1}) and (\ref{L4})  gives
\begin{equation}
\int\! {}_{2}F_{1}\left(\frac{a}{2},\frac{1}{2}-\frac{a}{2};1;\mathfrak{h}(x)\right)\,\mathrm{d}x=\frac{1}{2}x(1-x)(2-x){}_{2}F_{1}\left(\frac{a}{2}+1,\frac{3}{2}-\frac{a}{2};2;\mathfrak{h}(x)\right)+c.
\end{equation}
If we further define  $z=\mathfrak{h}(x)$ from which $x=1+\sqrt{1-z}$, the above integral becomes
\begin{equation}
\int\! (1-z)^{-\frac{1}{2}}{}_{2}F_{1}\left(\frac{a}{2},\frac{1}{2}-\frac{a}{2};1;z\right)\,\mathrm{d}z=z\sqrt{1-z}{}_{2}F_{1}\left(\frac{a}{2}+1,\frac{3}{2}-\frac{a}{2};2;z\right)+c
\end{equation}
which can be obtained as a special case of 1.15.3.9 in \cite{Prud}. A new indefinite integral involving a product of a Heun function with an incomplete elliptic integral of the first kind can be obtained by taking the function $h$ to be a solution of (\ref{2}) with $P$ given as in (\ref{pq}) and $\gamma=\delta=\epsilon=1/2$. Then, 3.131(3) in \cite{Grad} implies that
\begin{equation}
h(x)=\int_0^x\!\,\frac{\mathrm{d}u}{\sqrt{u(1-u)(a-u)}}=\frac{2}{\sqrt{a}}F\left(\varphi(x),\frac{1}{\sqrt{a}}\right),\quad
\varphi(x)=\arcsin{\sqrt{x}},
\end{equation}
where $F$ is the incomplete elliptic integral of the first kind. Furthermore, note that the above choice of the parameters $\gamma$, $\delta$ and $\epsilon$ requires that $\alpha+\beta=1/2$. Finally, we obtain the result
\[
\int\!\frac{\alpha(1-2\alpha)x-2q}{\sqrt{x(x-1)(x-a)}}F\left(\varphi(x),\frac{1}{\sqrt{a}}\right)H_\ell\left(a,q;\alpha,\frac{1}{2}-\alpha,\frac{1}{2},\frac{1}{2};x\right)\,\mathrm{d}x=
\]
\begin{equation}\label{Hn}
\sqrt{a}H_\ell\left(a,q;\alpha,\frac{1}{2}-\alpha,\frac{1}{2},\frac{1}{2};x\right)-
2\mathfrak{r}(x)F\left(\varphi(x),\frac{1}{\sqrt{a}}\right)H^{'}_\ell\left(a,q;\alpha,\frac{1}{2}-\alpha,\frac{1}{2},\frac{1}{2};x\right)+c.
\end{equation}
with $\mathfrak{r}(x)=\sqrt{x(x-1)(x-a)}$. 
If we let $z=\mathfrak{h}(x)$ with $x=1-\sqrt{1-z}$, make the identification $a=\alpha/2$, $b=1/4-a$, and $c=1/2$, and then, use (\ref{L1}) and (\ref{L4}), we obtain 
\[
\int\,z^{-\frac{1}{2}}(1-z)^{-\frac{1}{4}}F\left(\psi(z),\frac{1}{\sqrt{2}}\right){}_{2}F_{1}\left(a,\frac{1}{4}-a;\frac{1}{2};z\right)\,\mathrm{d}z=
\]
\begin{equation}\label{laut}
\frac{\sqrt{2}}{a(4a-1)}{}_{2}F_{1}\left(a,\frac{1}{4}-a;\frac{1}{2};z\right) +2z^{\frac{1}{2}}(1-z)^\frac{3}{4}F\left(\psi(z),\frac{1}{\sqrt{2}}\right){}_{2}F_{1}\left(a+1,\frac{5}{4}-a;\frac{3}{2};z\right)+c
\end{equation}
provided that $a\neq 0,1/4$ and $\psi(z)=\arcsin{\sqrt{1-\sqrt{1-z}}}$. To the best of our knowledge, the integral (\ref{laut}) seems to be new. We can also find indefinite integrals of products of hypergeometric functions with Heun functions by considering again (\ref{2}) with $\delta=0$ or $\epsilon=0$. The case $\delta$=0 yields $h(x)=x^{1-\gamma}{}_{2}F_{1}\left(\epsilon,1-\gamma;2-\gamma;x/a\right)$ and (\ref{1}) gives
\[
\int\!\frac{(x-a)^{\epsilon-1}(\alpha\beta x-q)}{x-1}{}_{2}F_{1}\left(\epsilon,\tau;1+\tau;\frac{x}{a}\right)H_\ell\left(a,q;\alpha,\beta,\gamma,0;x\right)\,\mathrm{d}x=(x-a)^\epsilon\cdot
\]
\[
\left\{\tau\left[{}_{2}F_{1}\left(\epsilon,\tau;1+\tau;\frac{x}{a}\right)+\frac{\epsilon x}{a(1+\tau)}{}_{2}F_{1}\left(\epsilon+1,1+\tau;2+-\tau;\frac{x}{a}\right)\right]H_\ell\left(a,q;\alpha,\beta,\gamma,0;x\right)\right.
\]
\begin{equation}\label{auf1}
\left.-x{}_{2}F_{1}\left(\epsilon,\tau;1+\tau;\frac{x}{a}\right)H^{'}_\ell\left(a,q;\alpha,\beta,\gamma,0;x\right)\right\}+c
\end{equation}
with $\epsilon=\alpha+\beta+\tau$ and $\tau=1-\gamma$. Finally, for $\epsilon$=0 we have $h(x)=x^{1-\gamma}{}_{2}F_{1}\left(\delta,1-\gamma;2-\gamma;x\right)$ and from (\ref{1}) we obtain
\[
\int\!\frac{(x-1)^{\delta-1}(\alpha\beta x-q)}{x-a}{}_{2}F_{1}\left(\delta,\tau;1+\tau;x\right)H_\ell\left(a,q;\alpha,\beta,\gamma,\delta;x\right)\,\mathrm{d}x=(x-1)^\delta\cdot
\]
\[
\left\{\tau\left[{}_{2}F_{1}\left(\delta,\tau;1+\tau;x\right)+\frac{\delta  x}{1+\tau}{}_{2}F_{1}\left(\delta+1,1+\tau;2+\tau;x\right)\right]H_\ell\left(a,q;\alpha,\beta,\gamma,\delta;x\right)\right.
\]
\begin{equation}\label{auf2}
\left.-x{}_{2}F_{1}\left(\delta,\tau;1+\tau;x\right)H^{'}_\ell\left(a,q;\alpha,\beta,\gamma,\delta;x\right)\right\}+c
\end{equation}
with $\delta=\alpha+\beta+\tau$ and $\tau=1-\gamma$. If we consider the special cases (\ref{L2}) and (\ref{L3}) for (\ref{auf1}) and (\ref{auf2})  together with relations (\ref{L5}) and (\ref{L6}), it can be shown after a lengthy but straightforward computation that the corresponding integrals are special cases of $1.15.2(4)$ in \cite{Prud}. This result can be interpreted as a further validation of formulae (\ref{auf1}) and (\ref{auf2}).\\
Furthermore, we can also take $h$ to be a solution of the ODE (\ref{3}) with $P$ and $Q$ as given in (\ref{pq}), i.e.
\begin{eqnarray}
h(x)&=&\mbox{exp}\left(-\int\!\frac{Q(x)}{P(x)}\,\mathrm{d}x\right)=\mbox{exp}\left(-\int\!\frac{\alpha\beta x-q}{K(x)}\,\mathrm{d}x\right),\quad
K(x)=\sum_{i=0}^2\mathfrak{k}_i x^i,\label{K}\\
\mathfrak{k}_2&=&\alpha+\beta+1,\quad
\mathfrak{k}_1=-[a(\gamma+\delta)+\alpha+\beta+1-\delta],\quad
\mathfrak{k}_0=a\gamma.
\end{eqnarray}
Let 
\begin{equation}\label{Delta}
\Delta=\mathfrak{k}_0\mathfrak{k}_2-\frac{1}{4}\mathfrak{k}_1^2.
\end{equation}
Then, the  integral in terms of which the function $h$ is expressed, can be explicitly computed by means of $2.103.5$ in \cite{Grad} yielding
\begin{equation}\label{ccc}
h(x)=\left\{
\begin{array}{ccc}
K^{-\frac{\alpha\beta}{2\mathfrak{k}_2}}(x)\mbox{exp}\left[\frac{\alpha\beta\mathfrak{k}_1+2q\mathfrak{k}_2}{2\mathfrak{k}_2\sqrt{\Delta}}\arctan{\left(\frac{2\mathfrak{k}_2 x+\mathfrak{k}_1}{2\sqrt{\Delta}}\right)}\right]&\mbox{if}~\Delta>0,\\
(x-x_0)^{-\frac{\alpha\beta}{\mathfrak{k}_2}}e^{\frac{C}{x-x_0}}  & \mbox{if}~\Delta=0,\\
K^{-\frac{\alpha\beta}{2\mathfrak{k}_2}}(x)\left[\frac{2\mathfrak{k}_2 x+\mathfrak{k}_1-2\sqrt{-\Delta}}{2\mathfrak{k}_2 x+\mathfrak{k}_1+2\sqrt{-\Delta}}\right]^{\frac{\alpha\beta\mathfrak{k}_1+2q\mathfrak{k}_2}{4\mathfrak{k}_2\sqrt{-\Delta}}}&\mbox{if}~\Delta<0,
\end{array}
\right.
\end{equation}
with $C=(\alpha\beta x_0-q)/\mathfrak{k}_2$. Note that we need to require that the singularity at $x_0=-\mathfrak{k}_1/2\mathfrak{k}_2$ lies outside the interval where the local solution of the Heun equation is defined. At this point (\ref{1}) leads to the following indefinite integral
\[
\int\!x^\gamma(x-1)^\delta(x-a)^\epsilon h(x) \frac{Q^2(x)+W(Q,P)(x)}{P^2(x)}H_\ell(a,q;\alpha,\beta,\gamma,\delta;x)\,\mathrm{d}x=
\]
\begin{equation}\label{hh1}
-x^\gamma(x-1)^\delta(x-a)^\epsilon h(x)\left[\frac{Q(x)}{P(x)}H_\ell(a,q;\alpha,\beta,\gamma,\delta;x)+H^{'}_\ell(a,q;\alpha,\beta,\gamma,\delta;x)\right]+c,
\end{equation}
where $P$ and $Q$ are given as in (\ref{pq}) and $W$ denotes the Wronskian. Equation (\ref{hh1}) gives rise to new indefinite integrals for the hypergeometric function. To this purpose let $a=2$, $q=\alpha\beta$, $\delta=\alpha+\beta-2\gamma+1$ which imply $\epsilon=\gamma$. Furthermore, if we define $z=2x-x^2$ from which $x=1-\sqrt{1-z}$, we get
\[
\int\!z^c(1-z)^{a+b-c}(1-\rho_1 z)(1-\rho_2 z)^{-2-\omega} {}_{2}F_{1}\left(a,b;c;z\right)\,\mathrm{d}z=
\]
\begin{equation}\label{hhh1}
\rho_3 z^c(1-z)^{a+b+1-c}(1-\rho_2 z)^{-\omega}\left[\frac{{}_{2}F_{1}\left(a,b;c;z\right)}{1-\rho_2 z}-{}_{2}F_{1}\left(a+1,b+1;c+1;z\right)\right]+\widetilde{c},
\end{equation}
where
\begin{eqnarray}
\rho_1&=&\frac{1+2(a+b+2ab)}{2[1-c+2(a+b+ab)]},\quad\rho_2=\frac{1+2(a+b)}{2c},\\
\rho_3&=&\frac{2c}{2(a+b+ab)-c+1},\quad\omega=\frac{2ab}{1+2(a+b)},
\end{eqnarray}
and $a=\alpha/2$, $b=\beta/2$, $c=\gamma$. Note that for $a=2$, $q=\alpha\beta$, $\delta=\alpha+\beta-2\gamma+1$ the integral (\ref{hhh1}) holds for both cases $\Delta>0$ and $\Delta<0$ because independently of the sign of the discriminant we have
\begin{equation}
h(x)=\left[(\alpha+\beta+1)(x^2-2x)+2\gamma\right]^{-\frac{\alpha\beta}{2(\alpha+\beta+1)}}.
\end{equation}
Let us consider the case $\Delta=0$ under the assumptions $a=2$, $q=\alpha\beta$, $\delta=\alpha+\beta-2\gamma+1$. Then,
\begin{equation}\label{lepanto}
\Delta=(\alpha+\beta+1)(\alpha+\beta+1-2\gamma).
\end{equation}
The equation $\Delta=0$ admits the following solutions
\begin{enumerate}
\item
$\gamma=(\alpha+\beta+1)/2$ if we look at $\Delta=0$ as an equation for $\gamma$.
\item
$\alpha_1=-\beta-1$ or $\alpha_2=2\gamma-\beta-1$ if we consider $\Delta=0$ as a quadratic equation in the parameter $\alpha$.
\item
$\beta_1=-\alpha-1$ or $\beta_2=2\gamma-\alpha-1$ if we consider $\Delta=0$ as a quadratic equation in the parameter $\beta$.
\end{enumerate}
Since $3.$ can be obtained from $2.$ by interchanging the parameters $\alpha$ and $\beta$, we will not consider this case. If $\gamma=(\alpha+\beta+1)/2$, then $\delta=0$ and
\begin{equation}
h(x)=(x-1)^{-\frac{\alpha\beta}{\alpha+\beta+1}}.
\end{equation}
Furthermore, we find that
\begin{equation}
\frac{Q^2(x)+W(Q,P)(x)}{P^2(x)}=\frac{\alpha\beta(\alpha\beta+\alpha+\beta+1)}{(\alpha+\beta+1)^2(x-1)^2},\quad
\frac{Q(x)}{P(x)}=\frac{\alpha\beta}{(\alpha+\beta+1)(x-1)}.
\end{equation}
Finally, if we define $z=2x-x^2$ from which $x=1-\sqrt{1-z}$ and make use of (\ref{L1}) and (\ref{L4}), we hand up with following known indefinite integral for the hypergeometric function
\[
\int\!z^{a+b+\frac{1}{2}}(1-z)^{-\frac{3}{2}-\xi}{}_{2}F_{1}\left(a,b;a+b+\frac{1}{2};z\right)\,\mathrm{d}z=\lambda z^{a+b+\frac{1}{2}}(1-z)^{-\xi}\cdot
\]
\begin{equation}\label{hhh1n}
\left[\frac{{}_{2}F_{1}\left(a,b;a+b+\frac{1}{2};z\right)}{\sqrt{1-z}}-\sqrt{1-z}{}_{2}F_{1}\left(a+1,b+1;a+b+\frac{3}{2};z\right)\right]+c,
\end{equation}
where
\begin{equation}
\lambda=\frac{2(2a+2b+1)}{2(2ab+a+b)+1},\quad\xi=\frac{2ab}{2a+2b+1},\quad a=\frac{\alpha}{2},\quad b=\frac{\beta}{2}.
\end{equation}
If we consider instead the case $\alpha=-\beta-1$, then $\delta=-2\gamma$ and $\epsilon=\gamma$. Moreover, we have
\begin{equation}
h(x)=\mbox{exp}\left(\frac{\beta(\beta+1)}{4\gamma}(x^2-2x)\right).
\end{equation}
Furthermore, we find that
\begin{equation}
\frac{Q^2(x)+W(Q,P)(x)}{P^2(x)}=\frac{\beta(\beta+1)\left[\beta(\beta+1)(x-1)^2+2\gamma\right]}{4\gamma^2},~
\frac{Q(x)}{P(x)}=\frac{\beta(\beta+1)(1-x)}{2\gamma}.
\end{equation}
Finally, if we define $z=2x-x^2$ from which $x=1-\sqrt{1-z}$ and make use of (\ref{L1}) and (\ref{L4}), we hand up with the following new indefinite integral for the hypergeometric function
\[
\int\!z^{c}(1-z)^{-c-\frac{1}{2}}(1-p_1 z)e^{-p_2 z}{}_{2}F_{1}\left(-b-\frac{1}{2},b;c;z\right)\,\mathrm{d}z=
\]
\begin{equation}\label{hhh1nT}
\lambda_1 z^c(1-z)^{\frac{1}{2}-c}e^{-p_2 z}\left[{}_{2}F_{1}\left(-b+\frac{1}{2},b+1;c+1;z\right)-{}_{2}F_{1}\left(-b-\frac{1}{2},b;c;z\right)\right]+\widetilde{c},
\end{equation}
where 
\begin{equation}
\lambda_1=\frac{2c}{b(2b+1)+c},\quad p_1=\frac{b(2b+1)}{b(2b+1)+c},\quad p_2=\frac{b(2b+1)}{2c},\quad b=\frac{\beta}{2},\quad c=\gamma.
\end{equation}
The case $\alpha=2\gamma-\beta-1$ will not be treated here because it gives rise to an integral similar to (\ref{hhh1n}). Further indefinite integrals can be obtained by choosing  the function $h$ to be a particular solution to (\ref{4}). There are several possibilities. For instance, if we make a s-homotopic transformation followed by a transformation of the independent variable in (\ref{4}), we find two linearly independent solutions of the form
\begin{equation}
h_i(x)=x^{-\alpha_i}(x-a)H_\ell\left(\frac{1}{a},q_i;\alpha_i,\beta_i,\gamma_i,0;\frac{1}{x}\right),\quad i=1,2
\end{equation}
with
\begin{eqnarray}
&&q_1=\frac{q-\alpha\beta}{a}+\omega-\rho,~
\alpha_1=\frac{1+\rho}{2},~\beta_1=\alpha_1+1,~\gamma_1=\frac{2(\omega-\alpha\beta+\rho)}{1+\rho},\\
&&q_2=\frac{\rho[\rho^2(q-\alpha\beta)-\alpha\beta(4a\omega+3)+3q]-4[\alpha^2\beta^2(a-3)+(3q-a+1)]}{a(1+\rho)^3},\\
&&\alpha_2=\frac{2\alpha\beta}{1+\rho},\quad
\beta_2=\alpha_2+1,\quad\gamma_2=2\alpha_2,\quad\rho=\sqrt{1-4\alpha\beta},\quad\omega=1-\alpha\beta\label{ro}
\end{eqnarray}
and (\ref{1}) gives rise to a couple of indefinite integral involving products of Heun functions and its derivatives, more precisely, for each $i=1,2$ we find
\[
\int\!x^{\gamma-1}(x-1)^{\delta-1}(x-a)^{\epsilon-1}K(x)h^{'}_i(x)
H_l(a,q;\alpha,\beta,\gamma,\delta;x)\,\mathrm{d}x=
\]
\begin{equation}\label{HeunHH}
x^{\gamma}(x-1)^{\delta}(x-a)^{\epsilon}
\left[h_i^{'}(x)H_l(a,q;\alpha,\beta,\gamma,\delta;x)-h_i(x)H^{'}_l(a,q;\alpha,\beta,\gamma,\delta;x)\right]+c,
\end{equation}
where $K$ has been defined in (\ref{K}). It is interesting to observe that if $q=0$ in (\ref{4}), the two linearly independent solutions are given by
\begin{equation}
h_1(x)=\frac{x-a}{(x-1)^{b_1}}{}_{2}F_{1}\left(a_1,b_1;c_1;\frac{a-1}{x-1}\right),\quad 
h_2(x)=\frac{x-a}{(x-1)^{a_1}}{}_{2}F_{1}\left(a_2,b_2;c_2;\frac{a-1}{x-1}\right),
\end{equation}
where
\begin{equation}
a_1=\frac{3-\rho}{2},\quad b_1=a_1-1,\quad c_1=2b_1,\quad
a_2=\frac{1+\rho}{2},\quad b_2=a_2+1,\quad c_2=2a_1
\end{equation}
with $\rho$ defined in (\ref{ro}), and from (\ref{HeunHH}) we can derive indefinite integrals of products of hypergeometric and Heun functions. Another possibility is to look at (\ref{4}) as a Heun equation with $\gamma=\delta=\epsilon=0$ and $\beta=-1-\alpha$. Then, we obtain an integral similar to (\ref{HeunHH}) but with $h_i^{'}(x)$ replaced by $H_l^{'}(a,q;\alpha,-1-\alpha,0,0;x)$. Furthermore, let us consider the particular solution $y(x)=H_l(a,q;\alpha,\beta,\gamma,\delta;x)$ to the Heun equation and a conjugate ODE to the Heun equation with $\overline{Q}(x)=(\alpha\beta x+q)/x(x-1)(x-a)$ for which we pick the particular solution $h(x)=H_l(a,-q;\alpha,\beta,\gamma,\delta;x)$. Then, (\ref{intcon}) yields the following indefinite integral involving products of Heun functions
\[
\int\!x^{\gamma-1}(x-1)^{\delta-1}(x-a)^{\epsilon-1}H_l(a,q;\alpha,\beta,\gamma,\delta;x)H_l(a,-q;\alpha,\beta,\gamma,\delta;x)\,\mathrm{d}x=
\]
\begin{equation}
\frac{x^{\gamma}(x-1)^{\delta}(x-a)^{\epsilon}}{2q}W\left(H_l(a,-q;\alpha,\beta,\gamma,\delta;x),H_l(a,q;\alpha,\beta,\gamma,\delta;x)\right)+c
\end{equation}
provided that $q\neq 0$. We conclude this section by constructing an indefinite integral involving products of Heun functions with complete elliptic integrals. To this purpose, we consider the Heun equation with $\gamma=1$ and $\delta=\epsilon=0$. Then, one particular solution is
\begin{eqnarray}
y(x)&=&\left( x-1 \right)^{\alpha}\psi(x),\label{y}\\
\psi(x)&=&H_l\left(1-a, \alpha^2(1-a)-\alpha-q,-\alpha,-\alpha+1,-2\,\alpha+1,0,
\frac{1-a}{1-x}\right).
\end{eqnarray}
As a conjugate equation we take
\begin{equation}
h^{''}(x)+\frac{1}{x}h^{'}(x)+\frac{1}{1-x^2}h(x)=0
\end{equation}
having a particular solution expressed in terms of complete elliptic integrals \cite{Grad}, namely $h(x)=\mathbf{E}(x^{'})$ with complementary modulus $x^{'}=\sqrt{1-x^2}$. Then, (\ref{intcon}) together with the functional relation 8.123(4) in \cite{Grad} between elliptic integrals
\begin{equation}
\frac{d\mathbf{E}(x^{'})}{dx^{'}}=\frac{\mathbf{E}(x^{'})-\mathbf{K}(x^{'})}{x^{'}}
\end{equation}
gives
\[
\int\!\frac{(x-1)^{\alpha-1}\mathfrak{Q}(x)}{(x+1)(x-a)}\mathbf{E}(x^{'})\psi(x)\,\mathrm{d}x=
\]
\begin{equation}
x(x-1)^\alpha\left[\frac{(1-\alpha)x-\alpha}{x^2-1}\mathbf{E}(x^{'})\psi(x)-\frac{x}{x^2-1}\mathbf{K}(x^{'})\psi(x)-\mathbf{E}(x^{'})\psi^{'}(x)\right]
\end{equation}
with $\mathfrak{Q}(x)=(1-\alpha^2)x^2-(a+q+\alpha^2)x-q$ and $\psi$ given by (\ref{y}).

\section{Comments and conclusions}
We applied the so-called Lagrangian method to obtain indefinite integrals of functions belonging to the family of Heun confluent functions. This approach allowed us to derive several novel indefinite integrals for the confluent, biconfluent, doubly confluent, and triconfluent Heun functions for which sample results have been provided. Our findings only scratch the surface of the wealth of new integral formulae one may obtain by using the aformentioned method.

\end{document}